\documentclass[11pt,a4paper]{amsart}
\usepackage[english]{babel}
\usepackage[latin1]{inputenc}
\usepackage{color}
\usepackage{graphics}

\usepackage{amsmath,amssymb,amstext,amsthm,amscd}

\def\ligne{\begin{center}
         \unitlength=1cm
         \begin{picture}(5,0.05)
         \line(1,0){5}
         \end{picture}
         \end{center}}

\newcounter{cimage}
\def\R{\mathbb R }

\def\G{\mathcal G }

\def\D{\cal D}
\def\grad{\mathrm{grad}}
\def\Eig{\mathrm{Eig}}

\def\k{\mathfrak k }
\def\g{\mathfrak g }

\def\vol{\mathrm{vol}}

\def\deg{\mathrm{deg}}

\def\tr{\mathrm{tr}}

\def\E{\mathrm{End}}
\def\rank{\mathrm{rang}}
\def\Herm{\mathrm{Herm}}

\def\Im{\mathrm{Im}}

\def\leq{\leqslant}

\newfam\msbfam

\font\tenmsb=msbm10
\font\sevenmsb=msbm10 at 7pt
\font\fivemsb=msbm10 at 5pt

\textfont\msbfam=\tenmsb
\scriptfont\msbfam=\sevenmsb
\scriptscriptfont\msbfam=\fivemsb

\def\textmap#1{\mathop{\vbox{\ialign{
                                 ##\crcr
     ${\scriptstyle\hfil\;\;#1\;\;\hfil}$\crcr
     \noalign{\kern-1pt\nointerlineskip}
     \rightarrowfill\crcr}}\;}}

\def\textlmap#1{\mathop{\vbox{\ialign{
                                 ##\crcr
     ${\scriptstyle\hfil\;\;#1\;\;\hfil}$\crcr
     \noalign{\kern-1pt\nointerlineskip}
     \leftarrowfill\crcr}}\;}}

\newfam\meuffam

\font\tenmeuf=eufm10
\font\sevenmeuf=eufm7

\textfont\meuffam=\tenmeuf
\scriptfont\meuffam=\tenmeuf
\scriptscriptfont\meuffam=\sevenmeuf

\def\germ{\fam\meuffam\tenmeuf}

\def\g{{\germ g}}

\def\EE{{\cal E}}
\def\FF{{\cal F}}

\def\tr{{\rm Tr}}

\def\Aut{{\rm Aut}}

\def\deg{{\rm deg}}
\def\Hom{{\rm Hom}}
\def\Aut{{\rm Aut}}

\def\Herm{{\rm Herm}}
\def\Vol{{\rm Vol}}

\def\id{{\rm id}}

\def\rk{{\rm rk}}

\def\cal{\mathcal}

\title{Optimal destabilizing  vectors in some  Gauge theoretical moduli problems}

\author{L. Bruasse}

\begin{document}
\newtheorem{toto}{toto}[section]
\def\thms{toto}


\theoremstyle{definition}
\newtheorem{definition}[\thms]{Definition}
\newtheorem{example}[\thms]{Example}
\newtheorem*{notations}{Notations}

\theoremstyle{plain}
\newtheorem{theorem}[\thms]{Theorem}
\newtheorem*{theorem2}{Theorem}
\newtheorem{proposition}[\thms]{Proposition}
\newtheorem{property}[\thms]{Property}
\newtheorem{lemma}[\thms]{Lemma}
\newtheorem{corollary}[\thms]{Corollary}

\theoremstyle{remark}
\newtheorem{remark}[\thms]{Remark}
\newtheorem*{notation}{Notation}
\theoremstyle{remark}
\newtheorem{exemple}[\thms]{Example}
\maketitle

\vspace{-0.7cm}
\begin{center}
 \small {\it IML, CNRS UPR 9016 \\ 163 avenue de Luminy,
 case 930, 13288 Marseille cedex 09 - France} \\
Email address: bruasse@iml.univ-mrs.fr
\end{center}


\ligne \vspace*{.2cm}
\begin{center}
\parbox{12cm}{%
\small \textsc{{Abstract.}}\ \ We show that  the ``classical'' Harder-Narasimhan filtration of a complex  vector bundle can be viewed as a limit object for the action of the gauge group, in the direction of an optimal destabilizing vector. This vector  appears as an extremal value of the so called  ``maximal weight function''. We give a complete description of these optimal destabilizing endomorphisms. Then we show that the same  principle holds for another important moduli problem: holomorphic pairs (i.e. holomorphic vector bundles coupled with morphisms with fixed source). We get a generalization of the Harder-Narasimhan filtration theorem for the associated notion of $\tau $-stability. These results suggest that the principle holds in the whole gauge theoretical framework.

\

}
\end{center}\vspace*{.2cm}

\ligne

\

\section{Introduction}
Harder and Narasimhan have proved \cite{hana} that any non semistable bundle on an algebraic curve admits a unique filtration by subsheaves such that the quotients are torsion free and semistable. This is now a classical result which was generalized for reflexive sheaves on projective manifolds \cite{sha},\cite{mar}, and then to any reflexive sheaf on an arbitrary compact Hermitian manifold \cite{bru1},\cite{bru2}.





The system of semistable quotients associated to  a non semistable vector bundle by the Harder-Narasimhan theorem  can be interpreted as a semistable object with respect to a new moduli problem~: the moduli problem for $G$-bundles, where $G$ is a product $\prod_i GL(r_i)$. 
The  motivation of this work is to find a general principle which applies to an arbitrary moduli problem: we want to associate in a {\emph{canonical way}} to a non semistable object  a new moduli problem and a semistable object for this new problem. One of the motivation is to find an analogous of the Harder-Narasimhan statement for other type of complex objects, for instance holomorphic bundles coupled with sections or endomorphisms.

\

The proper context for such a discussion is provided by the Geometric Invariant Theory (GIT) and its further developments (see \cite{kem}, \cite{mum}, \cite{ki}, \cite{ram}, \cite{t}, \cite{brute}). This approach gives indeed some powerful tools, such as the Hilbert criterium, to deal with the stability condition.
In \cite{brute}, we discuss the above  mentioned principle in the finite dimensional framework, for moduli problems associated with actions of a reductive group on a finite dimensional (\emph{possibly non compact}) manifold.  We look at a holomorphic action $\alpha :G\times F\to F$ of a complex reductive Lie group $G$ on a complex manifold $F$. We prove, that under a certain completeness assumption \cite{t}, the notion of stability  \emph{and of semistability } associated to the choice of a compact subgroup $K\subset G$ and of a $K$-equivariant  moment map $\mu $ may be defined in term  of a $G$-equivariant generalized maximal weight map $\lambda : H(G)\times F\to \R$, where $H(G)\subset \g$ is the union of all subspaces of the form $i\k$. In the classical GIT, this map is known as the Hilbert numerical function (see \cite{ki}).  We have proved the existence and unicity (up to equivalence) of an \emph{optimal destabilizing vector}  $s$ in the Lie algebra of $G$ associated to any non semistable point $f\in F$. The role of this  optimal vector corresponds precisely to that of  so-called adapted $1$-parameter subgroups in classical GIT (see \cite{kem} or \cite{ram}). We have shown that the path $t\mapsto e^{ts}f$ converges to a point $f_0$ which is semistable with respect to a natural action of the centralizator $Z(s)$ on submanifold of $F$. The assignment $f \mapsto f_0=\lim_{t\to \infty} e^{ts}f$ is the generalization of the Harder-Narasimhan statement.

Therefore the general principle may be stated as follows:  to get an analogue of the Harder-Narasimhan result for a given complex moduli problem, one has to give a gauge theoretical formulation of the problem (i.e. describe it in term of an action of a gauge group on certain complex variety), to compute the maximal weight map  and to study the optimal destabilizing vectors of the non semistable objects. Then the new semistable object is obtained as a limit in the direction of the optimal vector.

\

In gauge theory, there are well-known  links between the stability theory of vector bundles and the GIT concepts. The analogue of the maximal weight map in this setting has already been exploited by several authors, in particular to study  the Kobayashi-Hitchin correspondance (see for example  \cite{mut}, \cite{lute}). In the case of plain bundles, the relation between the algebraic geometric notion of stability and the Morse theory associated to the Yang-Mills functional (\cite{atibo},\cite{das}) is also well-known. {\emph{But what appears to be unavailable in the literature,  even in the case of plain bundles, is the very natural fact that the extremal value of the maximal weight map in the gauge theoretical framework corresponds exactly to the optimal direction which leads to the Harder-Narasimhan filtration  (or precisely to the associated system of semistable quotients)}}.     

\

The purpose of the present paper is to give a direct and complete proof of this fact in two classical  gauge theoretical problems : holomorphic vector bundles and holomorphic pairs (i.e. vector bundles coupled with morphisms with fixed source). These results were announced without proof in \cite{lute}. 

We first focus on the moduli problem of complex vector bundles: we prove that one can associate to any non semistable bundle a maximal destabilizing
element in the formal Lie algebra of the gauge group, then we give a complete description
of this optimal vector. The main tools here are  the notion of Harder-Narasimhan polygonal line (see \cite{bru}) and a well defined notion of energy for such lines.  The system of quotients of the ``classical'' Harder-Narasimhan
filtration appears as the  limit  in the direction of this optimal vector .

In a second part, we give an analogue description for the moduli problem associated to holomorphic pairs.
Here the corresponding notion of stability, called \emph{$\tau $-stability} (see \cite{bra}), depends on the choice of a real parameter. We prove, in this setting, a generalization of the Harder-Narasimhan filtration theorem for the $\tau$-stability (let us remark that a generalization of the Seshadri filtration in the context of $\tau$-stability was given a few years ago in \cite{bradaswen}).  Using an explicit formula for the maximal weight function, we describe the  optimal destabilizing vectors. Once again, the filtration is obtained as a limit in the direction of this extremal value of the maximal weight map. 

These two results suggest that the principle holds in the infinite dimensional gauge theoretical framework. It gives an intuitive way to build the ``Harder-Narasimhan object''  for a given moduli problem.

\section{Holomorphic vector bundles}\label{sec1}

We are interested here in the following gauge theoretical moduli problem: classifying the holomorphic structures on a given  complex vector bundle up to gauge equivalence. For clarity purposes, we will assume that the base manifold is a complex curve to avoid complications related to singular sheaves (and so the use of $L^2_1$-Sobolev sections).  Besides, the completeness property used in \cite{brute} is \emph{formally} satisfied for \emph{linear} moduli problems on a complex curve, so this is natural to work first in this framework.

\

 Let $E$ be a complex vector bundle of rank $r$ over the Hermitian curve  $(Y,g)$. We denote by $\G$ the complex gauge group $\G:=\Aut(E)$ whose formal Lie algebra is $A^0(\E(E))$.

Let $h$ be any Hermitian structure on $E$ and let us denote by 
$$ \cal K_h:=U(E,h) \subset \G$$
the real gauge group of unitary automorphisms of $E$ with respect to $h$.

We will use here the terminology developed in \cite{t} and \cite{brute}. An element $s\in A^0(\E(E))$ is said to be of {\emph{ Hermitian type}} if there exists a Hermitian metric $h$ on $E$, such that $s\in A^0(\Herm(E,h))$.

We identify a holomorphic structure $\EE$  on $E$ with the corresponding  integrable semiconnections $\overline{\partial}_\EE$  on $E$ (see \cite{lute}).
We are concerned with the stability theory for the action of $\G$ on the space $\cal H(E)$ of holomorphic structures. Fixing a Hermitian metric $h$, the moment map for the induced $\cal K_h$-action on $\cal H(E)$ is given by :
$$ \mu(  \EE)=\Lambda_g(F_{\EE,h}) + \frac{2\pi i}{\vol_g(Y)} m(\EE) \id_E$$
where $F_{\EE,h}$ is the curvature of the Chern connection associated to $\overline{\partial}_\EE$ and $h$ and 
$$m(\EE):=\frac{\deg (\EE)}{r}$$
is the slope of $\EE$.

Let us recall that  a holomorphic vector bundle $\EE\in \cal H(E)$ is semistable with respect to this moment map if and only if it is semistable in the sense of Mumford (\cite{mum}, \cite{lute}):

\

\emph{$m(\FF):=\frac{\deg(\FF)}{\rank(\FF)} \leqslant \frac{\deg(\EE)}{\rank(\EE)}:=m(\EE) \ \text{ for all reflexive subsheaves } \ \FF\subset \EE, \ \text{ such that } \ 0<\rank(\FF)<r.$}

\
   
One may give  an analytic  Hilbert type criterion for the stability theory associated to the moment map $\mu$. We will need the following notation : if $f$ is an endomorphism of a vector space $V$, we will put for any $a \in \R$
$$ V_f(a ):=\bigoplus_{a ' \leqslant a } \Eig(f,a ').$$
We extend the notation for endomorphisms of $E$ with constant eigenvalues in an obvious way.

Then, one has an explicit formula for the maximal weight map $\lambda $ (\cite{mu}) : if $\EE \in \cal H(E)$ and $s\in \Herm(E,h)$ then
$$\lambda^s({\cal E})=
\left\{\begin{array}{l}
\lambda_k\deg({\cal
E})+\sum\limits_{i=1}^{k-1}(\lambda_i-\lambda_{i+1})\deg({\cal
E}_i)-\frac{\deg({\cal E})}{r}\tr(s)\\ \hbox{if the eigenvalues}\
\lambda_1<\dots<\lambda_k
\
\hbox{of $s$ are}  \hbox{ constant and}\\ \hbox{$\EE_i:={\cal 
E}_s({\lambda_i})$ are holomorphic}\\ \\
\infty \ \hbox{if not}\ .
\end{array}\right.
$$

\

{\emph{{\bf Stability criterion:} A point $\EE \in \cal H(E)$ is semistable if and only if $\lambda^s(\EE)\geqslant 0$ for all $s\in A^0(\Herm(E,h))$}.

\

Let us come to the definition of    optimal destabilizing endomorphisms :

\begin{proposition}\label{promax}
Let $\EE$ be a non semistable bundle. There exists a Hermitian endomorphism $s_{op}\in A^0(\Herm(E,h)) $ such that
$$\lambda^{s_{op}}(\EE) = \inf_{\substack{s\in A^0(\Herm(E,h)) \\ \|s\|=1}} \lambda^s(\EE)$$
\end{proposition}

\begin{proof}

It is sufficient to consider the $s\in  A^0(\Herm (E,h))$ with  constant eigenvalues $\lambda _1<\cdots <\lambda _k$ and such that the $\EE_i$ are holomorphic. The condition $\|s\|=1$ implies that the eigenvalues $\lambda _i$ are bounded. Moreover, we know that the degree of a subbundle of $\EE$ is bounded above (\cite{bru1} prop 2.2), so that, writing $\lambda $ in the form  
$$\lambda ^s(\EE)= \deg(\EE) \big[\lambda _k-\frac{1}{r}\sum_{i=1}^k r_i \lambda _i\big] + \sum_{i=1}^{k-1} \deg(\EE_i) (\lambda _i - \lambda _{i+1}),$$
it becomes obvious that $\lambda^s(\EE) $ is bounded from below on the sphere $\|s\|=1$.

\

Now let $(s_n)_n$ be a sequence of Hermitian endomorphisms such that 
$$\lim_{n\to +\infty} \lambda ^{s_n}(\EE) = \inf_{\substack{s\in A^0(\Herm(E,h)) \\ \|s\|=1}} \lambda^s(\EE).$$
We always assume that the associated filtrations are holomorphic.  Going to a subsequence if necessary, we may suppose that each $s_n$ admits $k$  distinct eigenvalues $\lambda _1^n<\cdots <\lambda _k^n$. 

Let us recall the following result which is a direct consequence of the convergence theorem for subsheaves proved in \cite{bru2} (see also \cite{bru1} prop 2.9):
\begin{proposition} Let $(\mathcal F_n)_n$ be a sequence of subsheaves of $\mathcal{E}$ and assume that there exists a constant $c \in \R$ such that for all $n$ 
$$\deg(\mathcal F_n) \geqslant c$$
Then we may extract a subsequence $(\mathcal F_m)_m$ which converges \emph{in the sense of weakly holomorphic subbundles} to a subsheaf $\mathcal F$ of $E$. In particular $$\limsup_{m\to +\infty} \deg(\mathcal F_m) \leqslant \deg(\mathcal F).$$
\end{proposition}

Using this proposition and the  fact that the $\lambda _i^n$ are bounded, one can easily  extract from $(s_n)_n$ a subsequence $(s_m)_m$ such that :

\begin{enumerate}
{\it
\item there exist indices $0=j_0<j_1<\cdots <j_l=k$ and distinct values $\lambda _{j_1}< \cdots < \lambda _{j_l}$ such that for all $i\in \{j_p,\cdots,j_{p+1}\}$, $\lambda _i^m \to \lambda _{j_{p+1}}$;
\item there exists a filtration $0 =\EE_0 \subset \EE_{j_1}\subset \cdots \subset \EE_{j_l} =\EE$ such that each $\EE_{j_p}$ is the limit \emph{ in the sense of weakly holomorphic subbundles} of the $\EE_{j_p}^m$, which implies that $\deg (\EE_{j_p})\geqslant \limsup_{m\to +\infty} \deg(\EE_{j_p}^m)$.   
}
\end{enumerate}

For the second point, we use the fact that inclusion of subsheaves is preserved when going through the limit in the sense of weakly holomorphic subbundles \cite{bru1}.

Let $s$ be the Hermitian endomorphism whose eigenvalues are the $\lambda _{j_p}$ with corresponding filtration $\{\EE_{j_p}\}$ ($s$ may have less than $k$ distinct eigenvalues). Then we have
$$ \lambda ^s(\EE) \leqslant \liminf_{m\to +\infty}  \lambda^{s_m}(\EE) = \inf_{\substack{s\in A^0(\Herm(E,h)) \\ \|s\|=1}} \lambda^s(\EE)$$
so that $s$ is an optimal destabilizing endomorphism of $\EE$.
\end{proof}

Let us come to the proof of the result stated in \cite{brute} :

\begin{theorem}\label{the1}
Let $\EE \in \cal H(E)$ be a non semistable bundle. Then it admits a unique optimal destabilizing element $s_{op}\in A^0(\Herm(E,h))$ which is given by 
$$s_{op}=\frac{1}{\sqrt{\Vol(Y)}\sqrt{\sum\limits_{i=1}^k r_i\left[ \frac{\deg({\cal E}_{i}/{\cal
E}_{i-1})}{r_i }- \frac{\deg({\cal E})}{r}\right]^2 }}
\sum_{i=1}^k\left[\frac{\deg({\cal E})}{r}- \frac{\deg({\cal 
E}_{i}/{\cal E}_{i-1})}{r_i
}\right] \id_{F_i}\ ,
$$
where 
$$0=\EE_0\subset {\cal E}_1\subset {\cal E}_1\subset\dots\subset 
{\cal E}_k={\cal E}
$$
is the Harder-Narasimhan filtration of $\EE$, $F_i$ is the $h$-orthogonal complement of $\EE_{i-1}$ in $\EE_i$ and $r_i:=\rank(\EE_i/\EE_{i-1})$.
\end{theorem}
 
\begin{proof}
Let $s\in A^0(\Herm(E,h))$ such that $s$ has constant eigenvalues $\lambda_1< \lambda_2 <\cdots <\lambda _k$ and the filtration given by the $\EE_i=\EE_s(\lambda _i)$ is holomorphic. Let us denote by 
$$m_i:= \frac{\deg(\EE_i/\EE_{i-1})}{r_i}\ , 1\leqslant i \leqslant k, \  \ m:=m(\EE)$$
the slopes of the associated quotient sequence. The expression of the maximal weight map becomes :
$$ \lambda ^s(\EE)= \sum_{i=1}^k \lambda _i r_i (m_i-m).$$

We want to minimize this expression with respect to $s$ under the assumption 
$$\|s\|= \sqrt{\int_Y \tr(ss^\star).vol_g}=\sqrt{\Vol(Y)}\sqrt{\sum_{i=1}^k r_i \lambda _i^2}=1.$$
Assume first that the filtration $0=\EE_0 \subset \EE_1 \subset  \cdots \subset \EE_k =\EE$ is fixed and hence that  the $r_i$ and the $m_i$ are fixed. We have to minimize the map
$$g_{\{\EE_i\}}:(\lambda _1,\cdots, \lambda _k)\mapsto  \sum_{i=1}^k \lambda _i r_i (m_i-m)$$
over the smooth ellipsoid
$$ S=\big\{ (\lambda _1,\cdots,\lambda _k) \ | \  \sum_{i=1}^k r_i \lambda _i^2=1/\Vol(Y)\big\}.$$ 
with the additional open  condition 
$$ \lambda _1 < \lambda _2< \cdots <\lambda _k. \ \ (\star)$$
 
Resolving the Lagrange problem for $g|_S$ we see that there are two critical points of $g$ over $S$ which are obtained for $$\lambda _i= \varepsilon (m_i-m), \ 1\leqslant i\leqslant k \text{ and } \varepsilon =\pm 1 \frac{1}{\sqrt{\Vol(Y)}\sqrt{\sum_{i=1}^{k} r_i(m_i-m)^2}}.$$ 
Let us remark that $g_{\{\EE_i\}}$ is negative only for $\varepsilon <0$. 

We need the following
\begin{lemma}\label{gradlem}
The filtration $\{\EE_i\}_{1\leqslant i\leqslant k}$ associated to any optimal destabilizing element $s_{op}$  satisfies $m_1 > m_2 > \cdots > m_k$.
\end{lemma}

\begin{proof}
If the sequence $(m_i)_{ 1\leqslant i\leqslant k}$ is not decreasing, the critical point of $g_{\{\EE_i\}}$ which may correspond to a minimum does not satisfy the condition $(\star)$.  So let $s_{op}$ be an optimal destabilizing element and $\lambda _1< \cdots <\lambda _k$ its eigenvalues. Assume that for the corresponding filtration the sequence $(m_i)_{ 1\leqslant i\leqslant k}$ is not decreasing, then $\nu=(\lambda _1, \cdots ,\lambda _k)$ is not a critical point of $g=g_{\{\EE_i\}}$, so that the gradient $\grad_\nu(g|_S)$ is non zero. Therefore, moving slightly the point $\nu $ in the opposite direction, we may get a new point $\nu '\in S$ which still satisfies the open condition $(\star)$ and with $g(\nu ')<g(\nu )$. Thus, the corresponding Hermitian endomorphism $s'$ satisfies $\lambda ^{s'}(\EE)<\lambda ^{s_{op}}(\EE)$ which is a contradiction.   
\end{proof}

Keeping in mind this result, it is sufficient to consider endomorphisms whose associated filtration $\{\EE_i\}_{1\leqslant i\leqslant k}$  satisfies $m_1 > m_2>\cdots> m_k$. We will call such a filtration an \emph{admissible filtration}.
 
For admissible filtrations, the map $ g_{\{\EE_i\}}$ reaches its minimum for 
$$(\lambda _1,\cdots,\lambda_k)= \frac{1}{\sqrt{\Vol(Y)}\sqrt{\sum_{i=1}^{k} r_i(m_i-m)^2}}(m-m_1, \cdots ,m-m_k)$$
and 
$$ g_{\{\EE_i\}}(\lambda _1,\cdots,\lambda _k)=- \frac{1}{\sqrt{\Vol(Y)}}\sqrt{\sum_{i=1}^k r_i (m_i-m)^2}.$$

Thus, we  have to maximize $\sum_{i=1}^k r_i (m_i-m)^2$  among all admissible  filtrations.

\

Let us remind some fundamental property of the Harder-Narasimhan filtration (see \cite{bru}, \cite{bru1} for details):
\begin{proposition}
The Harder-Narasimhan filtration  is the unique filtration $0=\EE_0\subset \EE_1\subset \cdots \subset \EE_l=\EE$ such that 
\begin{enumerate}
\item each quotient $\EE_i/\EE_{i-1}$ is semistable for $1\leqslant i\leqslant l$;
\item the slope sequence satisfies
$$ m(\EE_i/\EE_{i-1}) < m (\EE_{i+1}/\EE_i) \ \text{ for } 1\leqslant i\leqslant l-1 .$$
\end{enumerate}
\end{proposition}

\

We may associate to any filtration $\cal D = (0=\EE_0\subset \EE_1\subset \cdots \EE_l = \EE)$ a polygonal line $\cal P(\cal D)$ in $\R^2$ defined by the sequence of points $p_i=(\rank \EE_i, \deg \EE_i)$ for $0\leqslant i\leqslant l$. The line associated to the Harder-Narasimhan filtration is called \emph{the Harder-Narasimhan polygonal line}.

\begin{proposition}[\cite{bru}]\label{brupol}
  For any subsheaf $\cal F$ of $\EE$, the point $(\rank \cal F,\deg \cal F)$ is located below  the Harder-Narasimhan polygonal line. As a consequence, any polygonal line associated to a filtration of $\EE$ is  located below the Harder-Narasimhan polygonal line. 
\end{proposition}

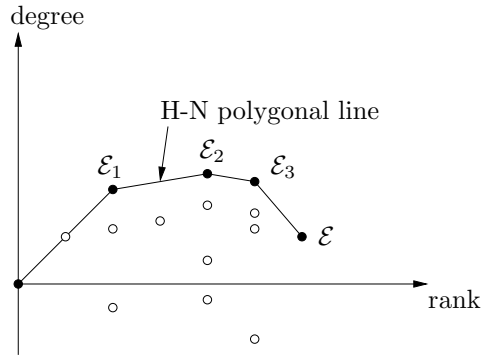
\centerline{\input{figure1.pstex_t}}

\

For any filtration $\cal D=\{\EE_i\}_{1\leqslant i\leqslant k}$, the condition $m_1 > \cdots > m_k$ is equivalent to the fact that the line  $\cal P(\cal D)$  is concave and let us remark that it is satisfied by the Harder-Narasimhan filtration. Now, the expression $\sum_{i=1}^k r_i (m_i-m)^2$ can be interpreted as an energy of the corresponding polygonal line. Indeed,  let us denote by $f_{\D}\in C^0([0,r],\R)$ the map whose graph is the polygonal line $\cal P(\D)$; this is a piecewise affine function and so  let $a_0=0<a_1 <\cdots <a_l=r$ be the corresponding partition of $[0,r]$.  We define
$$  E(f_\D):=\sum_{i=1}^k r_i (m_i-m)^2 = \sum_{i=0}^{l-1} \int_{a_i}^{a_{i+1}} (f'_\D (x) - m)^2 \mathrm{ d}x $$

\begin{proposition}
Let $\frak D$ be the Harder-Narasimhan filtration of $\EE$ and $\D$ be any other  admissible filtration,  then $E(f_\frak D) > E(f_\D)$. The maximum of the energy among all concave polygonal lines associated to a filtration of $\EE$ is obtained for the Harder-Narasimhan polygonal line.
\end{proposition}
 
\begin{proof}
Let $\frak D$ the Harder-Narasimhan filtration of $\EE$ and $\D$ any other filtration whose polygonal line is concave. 
This is essentially a result about concave piecewise affine maps :
\begin{lemma}\label{conaff}
Let $f,g :[0,r] \to \R$ be two continuous piecewise affine concave functions. Assume that $g(0)=f(0)$ and $f(r)=g(r)$. If  for all $x\in [0,r]$, $f(x)\geqslant g(x)$ then $E(f)\geqslant E(g)$. Equality may occur if and only if $f=g$. 
\end{lemma}

\begin{proof}

Let $a_0=0<a_1 <\cdots  <a_l=r$ be a partition of $[0,r]$ common to $f$ and $g$. Let $F(t,x)=t  f(x) + (1-t) g(x)$ then $x\to F(t,x)$ is continuous, concave,  and  affine on each segment $[a_i,a_{i+1}]$. So we may define for each $t$ the energy $E(F(t,x))$ of the map $x\mapsto F(t,x)$.

One has
\begin{align*}
\frac{d}{dt} E(F_t) = & \sum_{i=0}^{l-1}  2 \int_{a_i}^{a_{i+1}} \frac{d^2}{dtdx} F(t,x) (\frac{d}{dx} F(t,x) - m) \mathrm{d} x \\
= & \sum_{i=0}^{l-1}  [\frac{d}{dt} F(t,x) (\frac{d}{dx} F(t,x) - m)]_{a_i}^{a_{i+1}} \\ 
  & \ \ \ \ \ \ \ \ \ \ \ \ \ \ \ \  - 2 \sum_{i=0}^{l-1} \int_{a_i}^{a_{i+1}} \frac{d}{dt} F(t,x) \frac{d^2}{dx^2} F(t,x) \mathrm{d} x\\
 \geqslant &  \sum_{i=1}^{l-1} \big[ (f (a_i) - g (a_i))\times \\
 &  \ \ \ \ \ \ \  (t((f)'_g (a_i) -  (f)'_d (a_i)) + (1-t)((g)'_g (a_i) -  (g)'_d (a_i)))\big] \\
 > &  0
\end{align*}

Indeed, the concavity of each polygonal line implies that its derivative is
decreasing, that $\frac{d^2}{dx^2} F(t,x)\leqslant 0$ and we have of  course $(f(a_i) - g (a_i))\geqslant 0$.
 Let us remark that if $f\not= g$ the strict inequality $E(f) > E(g)$ obviously occurs.
\end{proof}

\begin{tabular}{cc}
\parbox{5cm}{%
\vspace*{-4cm}
Then, by proposition \ref{brupol} $\cal P(\D)$ lies below $\cal P(\frak D)$, that is for any $x\in [0,r]$, $f_\D (x)\leqslant  f_\frak D (x)$ with strict inequality on an open subset. We use the lemma to conclude.
}
& \input{figure2.pstex_t}
\end{tabular}

\end{proof}

Summarising the results, we have proved that $\lambda^s (\EE)$ reaches its minimum for a Hermitian element $s\in A^0(\Herm(E,h))$ whose associated filtration is the Harder-Narasimhan filtration of $\EE$ and whose eigenvalues are $\lambda _i= \frac{m-m_i}{\sqrt{\Vol(Y)}\sqrt{\sum_{i=1}^k r_i (m_i-m)^2}}$. Moreover it follows from the proof that this is a strict minimum.
\end{proof}

Now using the identification of $\cal H(E)$ with the space of integrable semiconnections, it is not difficult (see for instance \cite{mut} lemma 2.3.2) to show that $$\lim_{t\to +\infty} (e^{t{s_{op}}}). \overline{\partial}_\EE = \overline{\partial}_{F_1}\oplus \cdots \oplus \overline{\partial}_{F_k}.$$ 
In other words the holomorphic structure $e^{ts_{op}}\EE$ converges to the direct sum holomorphic structure $\bigoplus_{i=1}^k \EE_i/\EE_{i-1}$ as $t\to +\infty$. 

\emph{This illustrates our principle that the Harder-Narasimhan filtration is  a limit object precisely in the direction given by the extremal value of the maximal weight map}. The direct sum holomorphic structure is of course  semistable with respect to the smaller gauge group $\prod_{i=1}^k Aut(E_i/E_{i-1})$.

\section{Holomorphic pairs}
Here we will give an analogous result for a different  moduli problem associated to holomorphic pairs. 

Let $\cal F_0$ be a fixed holomorphic vector bundle of rank $r_0$ with a fixed Hermitian metric $h_0$ and $E$ a complex vector bundle of rank $r$ on the Hermitian curve $(Y,g)$. We are interested in the following moduli problem: classifying the holomorphic pairs $(\EE,\varphi )$ where $\EE$ is a holomorphic structure on $E$ and $\varphi $ is a holomorphic morphism $\varphi :\cal F_0 \to \EE$. Such a pair will be called a holomorphic pair of type $(E,\cal F_0)$ and we will denote by $\cal H(E,\cal F_0)$ the space of such pairs. Here the complex gauge group is once again  the group $\cal G:= \Aut(E)$.

Let us fix a Hermitian metric $h$ on $E$, and let us denote by $\cal K_h:=U(E,h)$ the group of unitary automorphisms.
The moment map for the $\cal K_h$ action on $\cal H(E,\cal F_0)$ has the form :
$$ \mu({\cal E},\varphi)=\Lambda_g F_{{\cal
E},h}-\frac{i}{2}\varphi\circ\varphi^*+\frac{i}{2} t\id_E\ .$$
In the sequel we will assume that the topological condition $\mu(\EE)\geqslant \tau$ holds.

Then we have the following characterisation of semistable pairs $(\EE,\varphi)$ (see \cite{bra}):
let $\tau :=\frac{1}{4\pi} t \Vol_g(Y)$, then $(\EE,\varphi)$ is semistable with respect to the moment map $\mu$ if and only if it is {\emph{ $\tau$-semistable}} in the following sense~:
\vspace{1.5mm}
\begin{enumerate}
\item
$m(\FF):=\frac{\deg(\FF)}{\rk(\FF)}\leq\tau$ for all reflexive subsheaves
$\FF\subset\EE$ with $0<\rk(\FF)<r$.
\vspace{2mm}
  \item $m(\EE/\FF):=\frac{\deg(\EE/\FF)}{\rk(\EE/\FF)}\geqslant \tau$
for all reflexive subsheaves
  $\FF\subset\EE$ with $0<\rk(\FF)<r$ and $\varphi\in H^0(\Hom(\FF_0,\FF))$.
\end{enumerate}

Now we may give an analogue of the Harder-Narasimhan theorem for this notion of stability:

\begin{theorem}\label{hardert}
Let $({\cal E},\varphi)$ be a non
$\tau$-semistable holomorphic pair of type $(E,{\cal F}_0)$. Then there exists a unique holomorphic filtration
with torsion free quotients
$$0=\EE_0\subset \EE_1
\subset\dots\subset\EE_m\subset\EE_{m+1}\subset\dots\subset\EE_k=\EE
$$
of $\EE$ such that:
\begin{enumerate}
  \item The
slopes sequence satisfies:
$$\hspace{-.6cm}m(\EE_{1}/\EE_0) >\dots
> m(\EE_{m}/\EE_{m-1})>\tau>
m(\EE_{m+2}/\EE_{m+1})>
\dots>m(\EE_{k}/\EE_{k-1})\ .$$
and the additional condition 
$$ \tau \geq  m(\EE_{m+1}/\EE_m)$$
\item The quotients $\EE_{i+1}/\EE_i$ are semistable  for $i\ne m$.
\item One of the following properties holds:
\begin{enumerate}
\item
$$\Im(\varphi)\not\subset{\cal E}_m\ ,\ \tau>
\frac{\deg(\EE_{m+1}/\EE_{m})}{\rk(\EE_{m+1}/\EE_{m })}$$
  and the pair
$(\EE_{m+1}/\EE_{m},\bar\varphi)$ is $\tau$-semistable, where $\bar\varphi$ is
the \\$\EE_{m+1}/\EE_{m}$-valued morphism induced by $\varphi$.
\item
$$\Im(\varphi)\not\subset{\cal E}_m\ ,\ \tau=
\frac{\deg(\EE_{m+1}/\EE_{m})}{\rk(\EE_{m+1}/\EE_{m })}
$$
and $\EE_{m+1}/\EE_{m}$ is semistable of slope $\tau$. This implies that the
pair $(\EE_{m+1}/\EE_{m},\bar\varphi)$ is $\tau$-semistable.
\item $\Im(\varphi)\subset {\cal E}_m$ and $\EE_{m+1}/\EE_{m}$ is semistable.
\end{enumerate}
\end{enumerate}

Moreover, in the cases  (b) and (c) the obtained filtration coincides with the
``classical'' Harder-Narasimhan filtration of ${\cal E}$ and the additional condition $m(\EE_{m+1}/\EE_m)>m(\EE_{m+2}/\EE_{m+1})$ holds.
\end{theorem}
\begin{proof}
We will use the results and methods of theorem 3.2 in \cite{bru1}. In order to build the filtration, let us first consider type (i) destabilizing subsheaves of $\EE$, i.e.  subsheaves $\FF$ which satisfy $m(\FF)>\tau$. If there exists such a subsheaf, we let  $\EE_1$ be the maximal destabilizing type (i) subsheaf: it is nothing but the first element of the Harder-Narasimhan filtration of $\EE$. If $\Im(\varphi )\subset \EE_1$ then we follow with the classical Harder-Narasimhan filtration (case (c)), if not we consider the pair $(\EE/\EE_1,\varphi_1)$ where $\varphi_1$ is induced by $\varphi $ and we follow the same principle until there is no more type (i) destabilizing subsheaves: we get a sequence 
$$0=\EE_0\subset \EE_1 \subset \cdots \subset \EE_m$$
which coincides with the first terms of the usual Harder-Narasimhan filtration and such that $\EE/\EE_m$ has no type (i) destabilizing subsheaf.
 
Then we consider type (ii) destabilizing subsheaves and we take the subsheaf $\FF_1$ containing $\Im(\varphi )$ and minimizing the slope $\mu(\EE/\FF_1)$ (we take $\FF_1$ of maximal rank with this property, it exists and it is unique by the arguments developed in \cite{bru1}). This will be the {\emph{last}} term of our filtration. By definition, we have $\mu(\EE/\FF_1)<\tau$.  Moreover $\FF_1$ contains $\EE_m$ by the following lemma.

\begin{lemma}
Let $\cal G\subset \EE$ be a maximal destabilizing subsheaf of type (i) and let $\FF\subset \EE$ be a type (ii) maximal destabilizing subsheaf, then $\cal G\subset \FF$.
\end{lemma}
\begin{proof}
Assume first that $\cal G\cap \FF=0$, then $(\FF+\cal G)/F \simeq \cal G$ and of course $\Im(\varphi )\subset \cal G + \FF$. Then we have 
$$\mu((\FF+\cal G)/\FF)=\mu(\cal G) >\tau >\mu(\EE/\FF)$$   
and using the following exact sequence
$$ 0 \xrightarrow{} (\FF+\cal G)/\FF \xrightarrow{} \EE/\FF \xrightarrow{} \EE/(\FF+\cal G) \xrightarrow{} 0$$
we get $\mu(\EE/(\FF+\cal G))<\mu(\EE/\FF)$ which is a contradiction.
Now assume that $\FF\cap \cal G$ is a non trivial subsheaf of $\cal G$ then using the following exact sequence
$$ 0 \xrightarrow{} \cal G\cap \FF \xrightarrow{} \cal G \xrightarrow{} \cal G/\cal G\cap \FF \xrightarrow{} 0$$
and the isomorphism $(\cal G+\FF)/\FF \simeq \cal G/\cal G\cap \FF$, we get 
$$\mu((\FF+\cal G)/\FF)=\mu(\cal G/\cal G\cap \FF) <\mu(\EE/\FF) <\tau <\mu (\cal G)$$
and $\mu(\cal G\cap \FF)>\mu (\cal G)$ which is a contradiction.
So we get $\cal G\cap \FF = \cal G$ so that $\cal G \subset \FF$.
\end{proof}
\begin{remark}
This lemma simply states that one can always build a Harder-Narasimhan filtration starting either from its first term or from its last term.
\end{remark}

Following the process, we get a sequence $\FF_l\subset \cdots \subset \FF_1\subset \EE$ where $\FF_l$ does contain $\EE_m$ and admits no type (ii) destabilizing subsheaf. It is quite easy to prove that each  quotient $\FF_i/\FF_{i+1}$ is  semistable.  Moreover, it is clear that $\FF_l/\EE_m$ has no type (i) destabilizing subsheaf. So, putting things together, we get a filtration
$$0\subset \EE_0\subset \cdots \subset \cdots \EE_m \subset \EE_{m+1}(=\FF_l) \subset \cdots \subset \EE_{k-1}=(\FF_1)\subset \EE_k=\EE$$
where $\Im(\varphi )\subset \EE_{m+1}$.

Also we have
$$\hspace{-.6cm}m(\EE_{1}/\EE_0) >\dots
> m(\EE_{m}/\EE_{m-1})>\tau$$
and
$$m(\EE_{m+2}/\EE_{m+1})>
\dots>m(\EE_{k}/\EE_{k-1})\ .$$

If $\Im(\varphi )\subset \EE_{m}$, then clearly the filtration coincides with the classical Harder-Narasimhan filtration and we are in  case (c) of the theorem. Else the pair $(\EE_{m+1}/\EE_m,\overline{\varphi})$ is by construction $\tau$-semistable with $m(\EE_{m+1}/\EE_m)\leqslant \tau$. Let us remark to conclude that if $m(\EE_{m+1}/\EE_m)= \tau$, then the notion of $\tau$-semistability and semistability coincide for $\EE_{m+1}/\EE_m$ (case (b) of the theorem).

Unicity is proved in the same way as in the algebraic ``classical'' case (see \cite{bru1} or \cite{bru}).
\end{proof}

\begin{remark}
Theorem \ref{hardert} clearly holds for any bundle over a compact Hermitian  manifold $(X,g)$ where $g$ is a Gauduchon metric (see \cite{bru1}). Indeed, we did not use the fact that $Y$ is one dimensional in this proof.
\end{remark}

Coming back to our gauge moduli problem, one can  once again give a formula for the maximal weight function:

$$\lambda^s({\cal E})=
\left\{\begin{array}{l} \lambda_k\deg({\cal
E})+\sum\limits_{i=1}^{k-1}(\lambda_i-\lambda_{i+1})\deg({\cal
E}_i)-\tau\tr(s)\\ \hbox{if the eigenvalues}\
\lambda_1<\dots<\lambda_k \ \hbox{of $s$ are}  \hbox{ constant,
$\EE_i:={\cal E}_s({\lambda_i)}$}\\ \hbox{are holomorphic, and
$\varphi\in
H^0(\Hom(\FF_0,\EE_s(0)))$}\ .\\
\\
\infty \ \hbox{if not}\ .
\end{array}\right.
$$

{\emph{ {\bf Criterion:} A holomorphic pair $(\EE,\varphi)$ is semistable with respect to the moment map $\mu$ if and only if $\lambda^s({\cal E})\geqslant 0$ for all $s\in A^0(\Herm(E,h))$.}}

\

Put again  $r_i:=\rk({\cal E}_{i}/{\cal E}_{i-1})$ and $m_i=m(\EE_i/\EE_{i-1})$.

We have the following result:
\begin{theorem}
For any
Hermitian metric $h$ on $E$, and any non semistable holomorphic pair $(\EE,\varphi)$, there exists a unique normalised Hermitian endomorphism $s_{op}\in A^0(\Herm(E,h))$ which satisfies
$$ \lambda^{s_{op}}(\EE)= \inf_{\substack{s\in A^0(\Herm(E,h)) \\ \|s\|=1}} \lambda^s(\EE).$$ 
It is given by
\begin{enumerate}
\item If $\Im(\varphi)\subset \EE_m$ then
$$s_{op}=\frac{1}{\sqrt{\Vol(Y)}\sqrt{\sum\limits_{i=1}^kr_i\left[
m_i- \tau\right]^2
}} \sum\limits_{i=1}^k\left[\tau- m_i\right]
\id_{F_i}\ .
$$
\item If $\Im(\varphi)\not\subset \EE_m$

$$s_{op}=\frac{1}{\sqrt{\Vol(Y)}\sqrt{\sum\limits_{\begin{array}{c}\scriptstyle
i=1\vspace{-1.5mm}\\
\scriptstyle i\ne m+1\end{array}}^kr_i\left[
m_i- \tau\right]^2
}} \sum\limits_{\begin{array}{c}\scriptstyle i=1\vspace{-1.5mm}\\
\scriptstyle i\ne m+1\end{array}}^k\left[\tau- m_i\right]
\id_{F_i}\ ,
$$
\end{enumerate}
where $F_i$ is the $h$-orthogonal complement of ${\cal E}_{i-1}$
in ${\cal E}_{i}$.
\end{theorem}

\begin{proof}
The existence of an optimal destabilizing element $s_{op}$ is proved in the same way as in proposition \ref{promax}: we simply use the fact that inclusion of subsheaves is preserved when going through the limit in the sense of weakly holomorphic subbundles to deal with  the additional condition $\Im(\varphi )\subset \EE_s(0)$.

\

Now, let $s\in A^0(\Herm(E,h))$ with constant eigenvalues $\lambda _1 <\cdots < \lambda _k$, such that the filtration given by the $\EE_i=\EE_s(\lambda _i)$ is holomorphic and $\Im(\varphi )\subset \EE_s(0)$. Using the previous notations the expression of the minimal weight map becomes 
$$\lambda ^s(\EE)=\sum_{i=1}^k \lambda _i r_i (m_i-\tau ).$$
As in the classical situation, we want to minimize this expression with respect to $s$ under the assumption
$$ \|s\|=\sqrt{\Vol(Y)}\sqrt{\sum_{i=1}^k r_i\lambda _i^2}=1$$
We use the same idea as in theorem \ref{the1}. Assume first that the filtration by eigenspaces is fixed;  we have to minimize the map
$$g_{\{\EE_i\}}:(\lambda _1,\cdots, \lambda _k)\mapsto  \sum_{i=1}^k \lambda _i r_i (m_i-\tau )$$
over the smooth ellipsoid
$$ S=\big\{ (\lambda _1,\cdots,\lambda _k) \ | \  \sum_{i=1}^k r_i \lambda _i^2=1/\Vol(Y)\big\}.$$ 
with the two  additional  conditions 
$$ \lambda _1 < \lambda _2< \cdots <\lambda _k. \  (\star 1) \ \ \text{ and }\ \ \Im(\varphi )\subset \EE_s(0) \  (\star 2).$$

Once again, resolving the Lagrange problem for ${g_{\{\EE_i\}}}|_S$, we see that there are two critical points of $g_{\{\EE_i\}}$ over $S$ which are obtained for 
$$\lambda _i= \varepsilon (m_i-\tau ), \ 1\leqslant i\leqslant k \text{ and } \varepsilon =\pm 1 \frac{1}{\sqrt{\Vol(Y)}\sqrt{\sum_{i=1}^{k} r_i(m_i-\tau )^2}}.$$ 
Let us remark that $g_{\{\EE_i\}}$ is negative only for $\varepsilon <0$. 

We have the following lemma:
\begin{lemma}
Assume $s$ is an optimal destabilizing Hermitian endomorphism and let  $\{\EE_i\}_{1\leqslant i\leqslant k}$  its associated filtration.
Then there exists $l$ such that 
\begin{enumerate}
\item $m_1>\cdots >m_l>\tau > m_{l+2}>\cdots >m_k=m(\EE)$;
\item $\tau \geq m_{l+1}$ and  the property $\Im(\varphi )\subset \EE_{l+1}$ holds.
\item if moreover $\Im(\varphi )\subset \EE_l$, then the additional condition $  m_{l+1}>m_{l+2}$ holds.
\end{enumerate}
A filtration which satisfies these conditions  will be called \emph{an admissible filtration}.
\end{lemma}

\begin{proof}
These are once again gradient arguments. Assume $s$ is optimal and let $l$ such that $\lambda_1< \cdots<\lambda_l<0 \leqslant \lambda_{l+1} <\cdots \lambda_k$.  

If $\Im(\varphi)\subset \EE_l$, since $\lambda_l<0$,  then the filtration must satisfy $m_1>\cdots >m_l>\tau \geqslant m_{l+1} > m_{l+2}>\cdots >m_k=m(\EE)$, otherwise the same gradient argument as in lemma \ref{gradlem} contradicts the optimality of $s$.
                  
If $\Im(\varphi)\not\subset \EE_l$, then obviously $ \Im(\varphi)\subset \EE_{l+1}$ and $\lambda_{l+1}=0$. Then restricting ${g_{\{\EE_i\}}}$ to $\{\lambda \, |\, \lambda_{l+1}=0\}$, a similar argument shows that the associated filtration satisfies i). For the second point, an explicit computation of the gradient of ${g_{\{\EE_i\}}}$, shows that if $m_{l+1}>\tau$, we may move the point $\lambda$ such that ${g_{\{\EE_i\}}}$ decreases and $\lambda_{l+1}\leqslant 0$, which contradicts the optimality of $s$. Thus $\tau \geq m_{l+1}$.
\end{proof}

So it is sufficient for our problem to minimize ${g_{\{\EE_i\}}}$ where $\{\EE_i\}$ is an admissible  filtration. Under this assumption  we get:
\begin{itemize}
\item If $\Im(\varphi) \subset \EE_l$, the minimal value of ${g_{\{\EE_i\}}}_{|_S}$ is obtained for 
$$(\lambda _1,\cdots,\lambda_k)= \frac{1}{\sqrt{\Vol(Y)}\sqrt{\sum_{i=1}^{k} r_i(m_i-\tau )^2}}(\tau -m_1, \cdots ,\tau -m_k)$$
and 
$$ g_{\{\EE_i\}}(\lambda _1,\cdots,\lambda _k)=- \frac{1}{\sqrt{\Vol(Y)}}\sqrt{\sum_{i=1}^k r_i (m_i-\tau )^2}.$$
\item If $\Im(\varphi )\not\subset \EE_l$, then  ${g_{\{\EE_i\}}}_{|_S}$ reaches its minimum for 
$$(\lambda _1,\cdots,\lambda_k)= \frac{1}{\|\lambda \|}\,(\tau -m_1, \cdots , \tau -m_{l}, 0, \tau -m_{l+2}, \cdots, \tau -m_k)$$
and 
$$ g_{\{\EE_i\}}(\lambda _1,\cdots,\lambda _k)=- \frac{1}{\sqrt{\Vol(Y)}}\sqrt{\sum_{\begin{array}{c}\scriptstyle i=1\vspace{-1.5mm}\\
\scriptstyle i\ne m+1\end{array}}^k r_i (m_i-\tau )^2}.$$
%
\end{itemize}

\

Now we want to minimize these expressions among all admissible filtrations. As in section \ref{sec1}, we will work on the polygonal line associated to any admissible filtration. We use the same notations as in section  \ref{sec1}, and we use the following definition for the energy $E(f)$ of a polygonal line $f$:
$$  E(f):=\sum_{i=1}^k r_i (m_i-\tau )^2 = \sum_{i=0}^{l-1} \int_{a_i}^{a_{i+1}} (f'_\D (x) - \tau )^2 \mathrm{ d}x $$

The difficult point here is that the polygonal line associated to an admissible filtration \emph{may no longer  be concave}.  We will in fact compare the energy of the different concave parts.

The following technical lemma will play an essential role in our proof:
\begin{lemma}\label{techlem}
Let $f:[0,r]\to \R$ and $g:[0,s]\to \R$ two distinct concave piecewise affine maps. Assume that 
\begin{enumerate}
\item $f(0)=g(0)$;
\item \label{cond2} $f(x)\geqslant g(x)$ for all $x\in [0,min(r,s)]$;
\item \label{cond3} $f'(x)\geqslant\tau $ and $g'(x)\geqslant\tau $ where they are defined;
\item $[f(r)-g(s)]/(r-s) \leqslant\tau $.
\end{enumerate}
Then $E(f)> E(g)$. 
\end{lemma}

\begin{proof}
Let us fix $f:[0,r]\to \R$; let $a_0=0< a_1<\cdots <a_l=s$ be the partition of $[0,s]$ corresponding to $g:[0,s]\to \R$. We will do an induction on the length $l$ of the partition.

For the case  $l=1$, assume first that $s\leqslant r$ then  the result is a consequence of lemma \ref{conaff} applied to $f|_{[0,s]}$ and $g$: using the additional condition \ref{cond3}, the conclusion of the lemma is still correct even if $g(s)\not=f(s)$. 

If $s>r$, let us denote by $h:[0,r]\to \R$ the straight line defined by $h(0)=0$ and $h(r)=f(r)$, then by lemma \ref{conaff}, $E(f)\geqslant E(h)$. Denote by $h_1$ the slope of the line $h$ and $g_1$ whose of the line $g$. By hypothesis  
$$(h(r)-g(s))/(r-s) \leqslant\tau $$
 such that 
$$h(r)-r\tau  \geqslant g(s)-s\tau $$
 and using conditions ii) and iii) we get 
$$\frac{1}{r}(h(r)-r\tau )^2 \geqslant \frac{1}{s}(g(s)-s\tau )^2,$$
so that 
$$r(h_1-\tau )^2\geqslant s(g_1-\tau )^2.$$
Thus $E(f)\geqslant E(h)\geqslant E(g)$, equality can only occur if $f=g$.

Now assume that  the result is proved for a  polygonal line of length less or equal than $l$ and let $g$ be a polygonal line of length $l+1$. Once again, if $s\leqslant r$, we use lemma \ref{conaff} to conclude. Else, we define a new polygonal line $h$ as follows: \begin{itemize}
\item $h_{|_{[0,a_l]}}=g_{|_{[0,a_l]}}$;
\item if adding the segment $[g(a_l),f(r)]$ preserves the concavity of $h$ we do so (case 1 below), if not we extend the last segment $[g(a_{l-1}),g(a_l)]$ in order that it reaches the line of slope $\tau $ going through $f(r)$ (case 2 below).
\end{itemize}

\

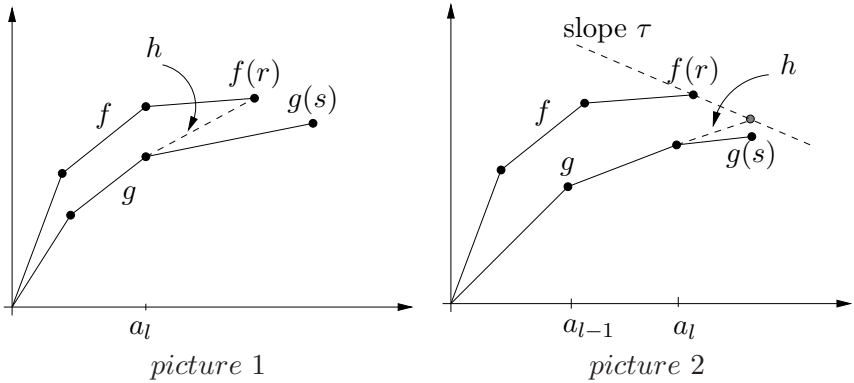
\centerline{$\begin{array}{cc}
\input{h1.pstex_t} & \input{h2.pstex_t}\\
picture \  1 & picture  \ 2
\end{array}$}

\

Using the same argument as in the first step of the induction we get easily $E(h)>E(g)$. Then we use lemma \ref{conaff} in the first case and the induction hypothesis applied to $f$ and $h$ in the second case to  get $E(f)\geqslant E(h)$.  
\end{proof}

 Let $\frak D=(0=\EE_0\subset \cdots\subset \EE_m\subset \EE_{m+1}\subset \cdots\subset \EE_k=\EE)$ be the generalized Harder-Narasimhan filtration given by theorem \ref{hardert} and $f_{\frak D}\in C^0([0,r],\R)$ the corresponding  piecewise affine map. Put $m_0=m$ if $\Im(\varphi )\subset \EE_m$, $m_0=m+1$ otherwise. Then $f_{\frak D}$ is not a concave map but admits two concave parts corresponding to $(0=\EE_0\subset \cdots\subset \EE_m)$ and $(\EE_{m_0}\subset \cdots\subset \EE_k=\EE)$. Keep in mind that the first $m$ subsheaves of this filtration  are just those of the classical Harder-Narasimhan filtration. Let $D=(0=\EE_0\subset \FF_1\subset \cdots \subset \FF_p=\EE)$ be any admissible filtration. Let us denote by $r'_i:=\rank(\FF_i/\FF_{i-1})$ and  $m'_i:=m(\FF_i/\FF_{i-1})$ and assume $m'_j>\tau \geqslant m'_{j+1}$. 

\

Assume that $\rank(\FF_j)>\rank(\EE_m)$. It follows by the definition of the filtration $\frak D$ that $(\deg(\EE_m)-\deg(\FF_j))/(\rank(\EE_m)-\rank(\FF_j))<\tau $: it is an obvious consequence of the fact that the point corresponding to $\FF_j$ is located below the classical Harder-Narasimhan polygonal line  and the fact that $\EE_m$ admits by definition no more type (i) destabilizing subsheaf. Hence applying  lemma \ref{techlem}  to the first concave part of each filtration, we get 
$$ \sum_{i=1}^j r'_i(\tau -m'_i)^2 \leqslant \sum_{i=1}^m r_i(\tau -m_i)^2.$$
If $\rank(\FF_j)\leqslant \rank(\EE_m)$, then the points corresponding to $\FF_i,1\leqslant i\leqslant j$ are located below the polygonal line $\cal P(\frak D)$ and using again lemma \ref{techlem} we get the same result.

\

We can apply exactly the same argument to the second part of the filtration. Let $j_0=j$ if $\Im(\varphi )\subset \FF_{j}$ and $j_0=j+1$ in the other case. Then, one can prove that the points corresponding to $(\FF_i)_{j_0\leqslant i\leqslant p}$ satisfy the following conditions :
\begin{itemize}
\item for any $i\in\{j_0,\cdots p\}$, if $\rank(\FF_i)>\rank(\EE_{m_0})$ then $\FF_i$ is located below the second part of the generalized Harder-Narasimhan polygonal line (i.e. the points corresponding to $\EE_{m_0}\subset \cdots \subset \EE$);
\item if $\rank(\FF_{j_0}) < \rank(\EE_{m_0})$, the relation 
$$(\deg(\FF_{j_0})-\deg(\EE_{m_0}))/(\rank(\FF_{j_0})-\rank(\EE_{m_0}))\geqslant \tau $$
holds;
\item the part of the associated polygonal line corresponding to the indices $\{j_0,\cdots,p\}$ is concave.
\end{itemize}
The first and the second point can be seen as a Harder-Narasimhan property for subsheaves containing $\Im(\varphi )$. Then using an analogue of lemma \ref{techlem} (in fact a symmetric proposition), we get:
$$ \sum_{u={j_0}}^p r'_u(\tau -m'_u)^2 \leqslant \sum_{u=m_0}^k r_u(\tau -m_u)^2.$$
This is exactly that we wanted since in the case where $\Im(\varphi)\not\subset \FF_j$, then the expression of $\lambda(\FF)$ \emph{does not} contain the $j^{th}$-term (the eigenvalue of the endomorphism must vanish).

In the following picture, the energy of the  black  and  grey parts of the filtrations are respectively compared to the energy of the corresponding parts of the generalized Harder-Narasimhan filtration. The segments drawn with  doted lines  give no contribution for the energy.  

\

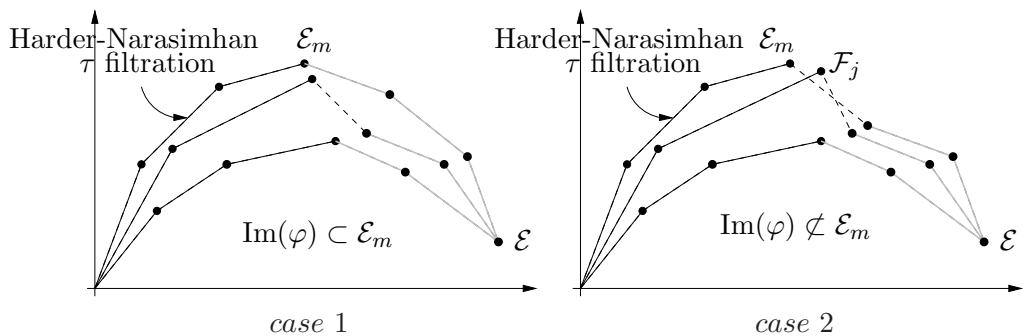
\centerline{$\begin{array}{cc}
\input{energy2.pstex_t} & \input{energy1.pstex_t}\\
case \  1 & case  \ 2
\end{array}$}

\

In either case, we have proved that the minimum of $\lambda (\EE)$ is achieved for the Hermitian element $s_{op}$ whose associated filtration is the generalized Harder-Narasimhan filtration with the corresponding eigenvalues described in the theorem.
\end{proof}

Once again, using the computation of \cite{mut}, we get:
\begin{itemize}
\item if $\Im(\varphi)\subset \EE_m$, then it is contained in the elements of the filtration corresponding to (strictly) negative eigenvalues of $s_{op}$, so that the path $t\mapsto e^{ts_{op}} (\EE,\varphi )$ converges as in the classical case  to the object
$$ (\EE_1/\EE_0,\cdots,\EE_k/\EE_{k-1}). $$
\item Else, if $\Im(\varphi )\not\subset \EE_m$, then it converges to the object
$$(\EE_1/\EE_0,\cdots,\EE_m/\EE_{m-1},(\EE_{m+1}/\EE_m,\overline{\varphi }),\EE_{m+2}/\EE_{m+1},\cdots,\EE_k/\EE_{k-1}).$$
where $\overline{\varphi }$ is the map induced from $\varphi $ (remark that in this case, the eigenvalue of $s_{op}$  corresponding to $F_m$ is vanishing). 
\end{itemize}

\

This illustrates once again our general principle.  By theorem \ref{hardert} the limit object is semistable with respect to the gauge group $\prod_{i=1}^k Aut(E_i/E_{i-1})$. 

This suggests that the principle holds in the whole (infinite dimensional) gauge theoretical  framework. 

\

\

\bibliographystyle{plain}
\bibliography{little}

\end{document}